\documentclass[a4paper]{article}
\usepackage{amsmath}
\usepackage{amssymb}
\usepackage{amscd}

\def\be{\begin{equation}}       \def\ee{\end{equation}}
\def\bd{\begin{displaymath}}    \def\ed{\end{displaymath}}
\def\beq{\begin{eqnarray}}      \def\eeq{\end{eqnarray}}
\def\bseq{\begin{eqnarray*}}    \def\eseq{\end{eqnarray*}}
\def\ba{\begin{array}}          \def\ea{\end{array}}
\def\ben{\begin{enumerate}}     \def\een{\end{enumerate}}

\def\lra{\longrightarrow}

\def\cop{\Delta}

\def\cnt{\varepsilon}
\def\ot{\otimes}

\def\lact{\triangleright}

\def\GLqtwo{{GL}_{q}(2)}

\def\GLpqtwo{{GL}_{p,q}(2)}
\def\GLhh'two{{GL}_{h,h'}(2)}

\def\Grss'{{G}_{r}^{s,s'}}
\def\Gmkk'{{G}_{m}^{k,k'}}

\def\rinv{{r}^{-1}}
\def\sinv{{s}^{-1}}

\def\ident{{\bf 1}}

\def\ca{{\mathcal A}}

\def\cd{{\mathcal D}}

\def\cl{{\mathcal L}}

\def\ct{{\mathcal T}}
\def\cu{{\mathcal U}}

\def\ta{\tilde{A}}
\def\tb{\tilde{B}}
\def\tc{\tilde{C}}
\def\td{\tilde{D}}
\def\tf{\tilde{F}}

\begin{document}
\begin{center}
{\Large \bf 
On the biparametric quantum deformation of $GL(2) \ot GL(1)$
}

\bigskip
\bigskip
\bigskip
\bigskip
\bigskip
\bigskip

{\large\bf Deepak Parashar
}

\smallskip
{\sl Max - Planck - Institute for Mathematics in the Sciences\\
Inselstrasse 22 - 26, D - 04103 Leipzig\\ 
Germany}\\

\smallskip
{\tt 
E-mail: Deepak.Parashar@mis.mpg.de
}

\bigskip
\bigskip
\bigskip
\bigskip
{\large \bf Abstract}
\end{center}
\medskip
We study the biparametric quantum deformation of $GL(2)\ot GL(1)$ and
exhibit its cross-product structure. We derive explictly the associated
dual algebra, i.e., the quantised universal enveloping algebra
employing the $R$-matrix procedure. This facilitates construction of a
bicovariant differential calculus which is also shown to have
a cross-product structure. Finally, a Jordanian analogue of the
deformation is presented as a cross-product algebra.
\par
\bigskip
\bigskip
\bigskip
\bigskip
\bigskip
\begin{center}
{\sl J. Math. Phys. {\bf 42} (2001) 5431 - 5443}
\end{center}
\newpage

\section{Introduction}

The biparametric quantum deformation of $GL(2) \ot GL(1)$ was introduced
in \cite{basu} as a novel Hopf algebra involving five generators $\{
a,b,c,d,f \}$ and two deformation parameters $\{ r,s \}$. From among the
five generators, four $\{ a,b,c,d \}$ correspond to $GL(2)$ and the fifth
one $f$ is related to $GL(1)$. These can be arranged in the matrix of
generators
\begin{equation}
\ct=\begin{pmatrix}f&0&0\\0&a&b\\0&c&d\end{pmatrix}
\end{equation}
with the labelling $0,1,2$. The associated solution of the quantum
Yang-Baxter equation is
\begin{equation}
R=\begin{pmatrix}
r & 0 & 0 & 0\\
0 & \mathsf{S}^{-1} & 0 & 0\\
0 & \Lambda & \mathsf{S} & 0\\
0 & 0 & 0 & R_{r}
\end{pmatrix}
\end{equation}
in block form, i.e., in the order $(00)$, $(01)$, $(02)$, $(10)$, $(20)$,
$(11)$, $(12)$, $(21)$, $(22)$ (which is chosen in conjunction with
the block form of the $\ct$-matrix) where
\[
R_{r}=\begin{pmatrix}
r & 0 & 0 & 0\\
0 & 1 & 0 & 0\\
0 & \lambda & 1 & 0\\
0 & 0 & 0 & r
\end{pmatrix}; \qquad
\mathsf{S}=\begin{pmatrix}
s & 0\\
0 & 1
\end{pmatrix}; \qquad
\Lambda=\begin{pmatrix}
\lambda & 0\\
0 & \lambda
\end{pmatrix}; \qquad
\lambda = r-r^{-1}
\]
The $R\ct\ct$ relations
\begin{equation}
R\ct_{1}\ct_{2}=\ct_{2}\ct_{1}R
\end{equation}
(where $\ct_{1}=\ct\ot \ident$ and $\ct_{2}=\ident \ot \ct$) give the
commutation relations between the generators $a$,$b$,$c$,$d$ and $f$
\begin{equation}
\begin{array}{ll}
ab=\rinv ba, \qquad & bd=\rinv db\\
ac=\rinv ca, \qquad & cd=\rinv dc\\
bc=cb, \qquad &[a,d]=(\rinv-r)bc
\end{array}
\label{grsabcd}
\end{equation}
and
\begin{equation}
\begin{array}{ll}
af=fa,&cf=sfc\\
bf=\sinv fb,&df=fd
\end{array}
\label{grsf}
\end{equation}
Note that the first set of these relations is exactly the $q$-deformation
of $GL(2)$ with deformation parameter $r$ while the second set involves
the fifth generator $f$ and the second deformation parameter $s$. This
results in a biparametric $q$-deformation of $GL(2)\ot GL(1)$, say,
$\ca_{r,s}$. The coproduct and counit is given as 
\begin{equation}
\begin{array}{l}
\cop (\ct) = \ct\dot{\ot} \ct \\
\cnt (\ct) = \ident
\end{array}
\end{equation}
The Casimir operator $\delta = ad-r^{-1}bc$ is invertible and determines
the antipode
\begin{equation}
S(f) = f^{-1}, \quad
S(a) = \delta^{-1}d, \quad
S(b) = -\delta^{-1}rb, \quad
S(c) = -\delta^{-1}r^{-1}c, \quad
S(d) = \delta^{-1}a
\end{equation}
The quantum determinant $\cd=\delta f$ is group-like but not central. Some
of the interesting features of the above quantum deformation are the
following:
\begin{itemize}

\item
If we write the set of generators $\{ a,b,c,d,f \}$ as $\{
f^{N}a,f^{N}b,f^{N}c,f^{N}d \}$ ($N$ being a fixed nonzero
integer), i.e., reducing the five-dimensional set to the
four-dimensional set, then we obtain an exact realisation of the
biparametric $(p,q)$-deformation of $GL(2)$, i.e., $\GLpqtwo$ subject
to the relations
\begin{equation}
p=r^{-1}s^{N} \quad \text{and} \quad q=r^{-1}s^{-N} 
\end{equation}
This realisation also reproduces the full Hopf algebraic structure
underlying $\GLpqtwo$.

\item
Another interesting feature of the $\ca_{r,s}$ deformation is that it can
be contracted (by means of the contraction procedure \cite{agha} based on
the concept of singular limit of a similarity transformation) to yield the
corresponding biparametric Jordanian deformation of $GL(2)\ot GL(1)$, 
which
in turn provides a complete realisation of the biparametric
$(h,h')$-deformation of $GL(2)$, i.e., $\GLhh'two$ in a manner similar
to that for the $q$-deformed case \cite{deeps}.

\item
Both the biparametric quantum and Jordanian deformations of $GL(2)\ot
GL(1)$ admit coloured extensions \cite{deeps} which also commute with the
contraction procedure.

\item
The physical interest in studying $\ca_{r,s}$ lies in the observation that
when endowed with a $\ast$-structure, this specialises to its compact
form, i.e., provides a biparametric $q$-deformation of $SU(2)\ot U(1)$
which is precisely the gauge group for the theory of electroweak
interactions.

\end{itemize}

Another deformation similar to $\ca_{r,s}$ has also been recently given in
\cite{grosse}, though in a different context. In the present article, we
give an explicit description of the algebra dual to $\ca_{r,s}$ as a
starting point in further investigation of this quantum group 
structure. Motivated by the relation of this deformation with gauge
theory, we also construct a bicovariant differential calculus since gauge
theories have an obvious differential geometric description. This would
then provide insights into possible scenarios for constructing $q$-gauge
theories based on this deformation. In pursuing our aim, we follow the
convenient $R$-matrix approach \cite{frt,jurco}. In Sec. II, we give the
cross-product structure and go over to the $R$-matrix duality in
Sec. III. The constructive calculus is presented in Sec. IV, while Sec. V
is a brief description of the Jordanian analogue. The results are
discussed in Sec. VI.

\section{Cross-product structure}

The biparametric $q$-deformation $\ca_{r,s}$ can also be
considered as the semidirect or cross-product 
$GL_{r}(2)\underset{s}{\rtimes} \mathbb{C}[f,f^{-1}]$
built on the vector space $GL_{r}(2)\ot \mathbb{C}[f,f^{-1}]$
where $GL_{r}(2)=\mathbb{C}[a,b,c,d]$ modulo the relations 
(\ref{grsabcd}) and $\mathbb{C}[f,f^{-1}]$ has the cross relations
(\ref{grsf}). Then, $\ca_{r,s}$ can also be interpreted as a skew Laurent
polynomial ring $GL_{r}[f,f^{-1};\sigma]$ where $\sigma$
is the automorphism given by the action of element $f$ on
$GL_{r}(2)$. Knowing properties of cross-product algebras
(general theory given in \cite{majid:book, ks}), we already know that the
algebra dual to $\ca_{r,s}$ would be the cross-coproduct coalgebra 
$\cu_{r,s}=U_{r}(gl(2))\underset{s}{\rtimes} \mathbb{C}[[\phi]]$ with
$\phi$ as an element dual to $f$. If we let $A=GL_{r}(2)$ and
$H=\mathbb{C}[f,f^{-1}]$, then $A$ is a left $H$-module algebra and the
action of $f$ on $GL_{r}(2)$ is given by
\begin{equation}
f\lact a=a,\qquad f\lact b=sb,\qquad f\lact c=s^{-1}c,\qquad f\lact d=d
\end{equation}
As a vector space, the dual is $\cu_{r,s}=U_{r}(gl(2))\ot U(u(1))$. Now,
the duality relation between $\langle GL_{r}(2), U_{r}(gl(2))\rangle$ is
already well-known \cite{sud}, while that between $\langle
\mathbb{C}[f,f^{-1}], U(u(1))\rangle$ is given by $\langle f, \phi \rangle
= 1$,
i.e., $U(u(1))=\mathbb{C}[[\phi]]$. More precisely, we work algebraically
with $\mathbb{C}[s^{\phi},s^{-\phi}]$ where $\langle f, s^{\phi} \rangle =
s$ (this is a standard notational convention which we adopt). This induces
duality on the vector space tensor products and the left action dualises
to the left coaction. This results in the dual algebra being a
cross-coproduct $\cu_{r,s}=U_{r}(gl(2))\underset{s}{\bf \rtimes}
\mathbb{C}[[\phi]]$. Let us recall \cite{sud} that $U_{r}(gl(2))$, the
algebra dual to $GL_{r}(2)$,
is isomorphic to the tensor product $U_{r}(sl(2))\ot \tilde{U}(u(1))$
where $U_{r}(sl(2))$ has the usual generators $\{ H,X_{\pm}\}$ and
$\tilde{U}(u(1))=\mathbb{C}[[\xi]]=\mathbb{C}[r^{\xi},r^{-\xi}]$ with
$\xi$ central. Therefore, $\cu_{r,s}$ is nothing but $U_{r}(sl(2))$ and
two central generators $\xi$ and $\phi$, where $\xi$ is the
generating element of $\tilde{U}(u(1))$ and $\phi$ is the
generating element of $U(u(1))$. Also note that $s^{\phi}$ ($s$ being 
the second deformation parameter) is dually paired with the
element $f$ of $\ca_{r,s}$. Defining the left coaction $U_{r}(gl(2)) \lra
U(u(1)) \ot U_{r}(gl(2))$ we have
\begin{equation}
X_{+}\lra s^{\phi}\ot X_{+}, \qquad X_{-}\lra s^{-\phi}\ot X_{-}, \qquad
H\lra 1\ot H, \qquad \xi\lra 1\ot \xi
\end{equation}
It can be checked that this gives the correct duality pairings. For
example, we have for $X_{+}$
\begin{equation}
\begin{array}{ll}
&\langle\Delta_{L}(X_{+}), 1\ot
\left( \begin{smallmatrix}a&b\\c&d\end{smallmatrix}\right) \rangle=
\langle s^{\phi}\ot X_{+}, 1\ot 
\left( \begin{smallmatrix}a&b\\c&d\end{smallmatrix}\right) \rangle=
\langle s^{\phi}, 1\rangle \langle X_{+},
\left( \begin{smallmatrix}a&b\\c&d\end{smallmatrix}\right) \rangle=
\left( \begin{smallmatrix}0&1\\0&0\end{smallmatrix}\right) \\
&\langle\Delta_{L}(X_{+}), f\ot 
\left( \begin{smallmatrix}a&b\\c&d\end{smallmatrix}\right) \rangle=
\langle s^{\phi}\ot X_{+}, f\ot
\left( \begin{smallmatrix}a&b\\c&d\end{smallmatrix}\right) \rangle=
\langle s^{\phi}, f\rangle \langle X_{+},     
\left( \begin{smallmatrix}a&b\\c&d\end{smallmatrix}\right) \rangle=
s \left( \begin{smallmatrix}0&1\\0&0\end{smallmatrix}\right) \\
&\langle X_{+}, f\lact
\left( \begin{smallmatrix}a&b\\c&d\end{smallmatrix}\right) \rangle=
\langle X_{+},
\left( \begin{smallmatrix}a&sb\\s^{-1}c&d\end{smallmatrix}\right) \rangle=
s \left( \begin{smallmatrix}0&1\\0&0\end{smallmatrix}\right)
\end{array}
\end{equation}
Therefore, the coalgebra structure of $\cu_{r,s}$ is given as
\begin{eqnarray}
& &\cop(X_{+}) = X_{+}\ot r^{\frac{H}{2}}+r^{-\frac{H}{2}}s^{\phi}\ot
X_{+}\\
& &\cop(X_{-}) = X_{-}\ot r^{-\frac{H}{2}}+r^{\frac{H}{2}}s^{-\phi}\ot
X_{+}\\
& &\cop(H) = H\ot 1+1\ot H\\
& &\cop(\xi) = \xi\ot 1+1\ot \xi\\
& &\cop(\phi) = \phi\ot 1+1\ot \phi
\end{eqnarray}
In this way, we have obtained the Drinfeld-Jimbo form of the dual algebra
$\cu_{r,s}$ using the cross-product construction. Given other approaches
to the problem of duality for quantum groups, we also construct explicitly
the dual algebra using the $R$-matrix procedure.

\section{$R$-matrix duality}

The biparametric $(r,s)$-deformation, $\ca_{r,s}$, of $GL(2)\ot GL(1)$ has
been defined in the previous section at the group level, i.e., as the
$q$-deformation of algebra of functions on $GL(2)\ot GL(1)$. In this
section, we derive explicitly the corresponding quantised universal
enveloping algebra, i.e., its dual within the framework of the
$R$-matrix formulation. We first construct functionals
(matrices) $\cl^{+}$ and $\cl^{-}$ which are dual to the matrix of
generators in the fundamental representation.  The linear functionals
$(\cl^{\pm})^a_{b}$ (following the method of \cite{frt, majid:book}) are
defined by their value on the elements of the matrix of generators $\ct$
\begin{equation}
\langle (\cl^{\pm})^a_{b}, \ct^c_{d} \rangle = (R^{\pm})^{ac}_{bd}
\label{lplusminus}
\end{equation}
where
\begin{eqnarray}
(R^{+})^{ac}_{bd} &=& c^{+} (R)^{ca}_{db}\\
(R^{-})^{ac}_{bd} &=& c^{-} (R^{-1})^{ac}_{bd}
\end{eqnarray}
and $c^{+}$ , $c^{-}$ are free parameters. Matrices $(\cl^{\pm})^a_{b}$
satisfy
\begin{equation}
\langle (\cl^{\pm})^a_{b}, uv \rangle =
\langle (\cl^{\pm})^a_{c}\ot (\cl^{\pm})^c_{d}, u\ot v \rangle =
(\cl^{\pm})^a_{c}(u) (\cl^{\pm})^c_{d}(v)
\label{functionals}
\end{equation}
\[
\text{i.e.} \quad
\cop (\cl^{\pm})^a_{b} = (\cl^{\pm})^a_{c}\ot (\cl^{\pm})^c_{b}
\]
For $\ca_{r,s}$, the $(R^{+})$ and $(R^{-})$ matrices read
\begin{equation}
(R^{+})= c^{+} \begin{pmatrix}
r & 0 & 0 & 0\\
0 & \mathsf{S} & \Lambda & 0\\
0 & 0 & \mathsf{S}^{-1} & 0\\
0 & 0 & 0 & R^{T}_{r}
\end{pmatrix}; \quad
(R^{-})= c^{-} \begin{pmatrix}
r^{-1} & 0 & 0 & 0\\
0 & \mathsf{S} & 0 & 0\\
0 & -\Lambda & \mathsf{S}^{-1} & 0\\
0 & 0 & 0 & R^{-1}_{r}
\end{pmatrix}
\end{equation}
where $R_{r}$, $\Lambda$ and $\mathsf{S}$ are the same as before and 
$R^{-1}_{r}=R_{r^{-1}}$. Before proceeding further, it is pertinent to
make the following remark about the $\cl^{\pm}$ functionals. Let $A(R)$ be
a bialgebra or a Hopf algebra underlying a $3\times 3$ quantum matrix and
let $\tilde{U}(R)$ be a similar matrix bialgebra with two full matrices
$\cl^{\pm}$ of generators. These may be
viewed as functionals $A(R)\lra \mathbb{C}$ via (\ref{lplusminus}), but
duality pairing at this level may be degenerate. So, we look at
appropriate quotients of these such that the pairing is non-degenerate. In
our case, upon quotienting $A(R)$ would descend to $\ca_{r,s}$,
and likewise $\tilde{U}(R)$ to the dual of $\ca_{r,s}$. The quotient on
$A(R)$ is obtained by setting certain entries of the T-matrix to zero. The
most general $3\times 3$ quantum matrix has nine elements
\begin{equation}
\ct=\begin{pmatrix}
\ct^{0}_{0} & \ct^{0}_{1} & \ct^{0}_{2}\\
\ct^{1}_{0} & \ct^{1}_{1} & \ct^{1}_{2}\\
\ct^{2}_{0} & \ct^{2}_{1} & \ct^{2}_{2}
\end{pmatrix}
\end{equation}
Now, let $\ct^{0}_{1} = 0 = \ct^{0}_{2}$ and $\ct^{1}_{0} = 0 =
\ct^{2}_{0}$. Checking the coideal property (via coproduct of $\ct$), we
have
\begin{equation}
\begin{array}{l}
\cop(\ct^{0}_{1})=\ct^{0}_{0}\ot \ct^{0}_{1} + \ct^{0}_{1}\ot \ct^{1}_{1}
+ \ct^{0}_{2}\ot \ct^{2}_{1}\\
\cop(\ct^{0}_{2})=\ct^{0}_{0}\ot \ct^{0}_{2} + \ct^{0}_{1}\ot \ct^{1}_{2}  
+ \ct^{0}_{2}\ot \ct^{2}_{2}\\
\cop(\ct^{1}_{0})=\ct^{1}_{0}\ot \ct^{0}_{0} + \ct^{1}_{1}\ot \ct^{1}_{0}
+ \ct^{1}_{2}\ot \ct^{2}_{0}\\
\cop(\ct^{2}_{0})=\ct^{2}_{0}\ot \ct^{0}_{0} + \ct^{2}_{1}\ot \ct^{1}_{0} 
+ \ct^{2}_{2}\ot \ct^{2}_{0}
\end{array}
\end{equation}
These generate biideals. Therefore, setting them to zero gives the
quotient of $A(R)$
\begin{equation}
\ct=\begin{pmatrix}
\ct^{0}_{0} & 0 & 0\\
0 & \ct^{1}_{1} & \ct^{1}_{2}\\
0 & \ct^{2}_{1} & \ct^{2}_{2}
\end{pmatrix}
= \begin{pmatrix}
f & 0 & 0\\          
0 & a & b\\
0 & c & d
\end{pmatrix}
= \ct (\ca_{r,s})
\end{equation}
Similarly, the quotient on $\tilde{U}(R)$ is obtained by setting certain
entries of $\cl^{\pm}$ matrices to zero. Starting with
\begin{equation}
\cl^{+}=\begin{pmatrix}
\cl^{+0}_{\enspace 0} & \cl^{+0}_{\enspace 1} & \cl^{+0}_{\enspace 2}\\
\cl^{+1}_{\enspace 0} & \cl^{+1}_{\enspace 1} & \cl^{+1}_{\enspace 2}\\
\cl^{+2}_{\enspace 0} & \cl^{+2}_{\enspace 1} & \cl^{+2}_{\enspace 2}
\end{pmatrix}, \qquad
\cl^{-}=\begin{pmatrix}
\cl^{-0}_{\enspace 0} & \cl^{-0}_{\enspace 1} & \cl^{-0}_{\enspace 2}\\
\cl^{-1}_{\enspace 0} & \cl^{-1}_{\enspace 1} & \cl^{-1}_{\enspace 2}\\
\cl^{-2}_{\enspace 0} & \cl^{-2}_{\enspace 1} & \cl^{-2}_{\enspace 2}
\end{pmatrix}
\end{equation}
we make the ansatz
\begin{equation}
\begin{array}{ll}
& \cl^{+2}_{\enspace 1} = 0 = \cl^{-1}_{\enspace 2}\\
& \cl^{+0}_{\enspace 1} = \cl^{+0}_{\enspace 2} = \cl^{+1}_{\enspace 0} =
\cl^{+2}_{\enspace 0} = 0\\
& \cl^{-0}_{\enspace 1} = \cl^{-0}_{\enspace 2} = \cl^{-1}_{\enspace 0} =   
\cl^{-2}_{\enspace 0} = 0
\end{array}
\end{equation}
and, similar to the above for $A(R)$, check the coideal property. We
also verify explicitly \cite{majid:book} that this ansatz is compatible
with the
duality pairing
\begin{equation}
\begin{array}{ll}
& \langle \cl^{+2}_{\enspace 1}, \ct^{i}_{j}\rangle= R^{+2i}_{\enspace 1j}
=R^{i2}_{j1} = 0\\
& \langle \cl^{-1}_{\enspace 2}, \ct^{i}_{j}\rangle= R^{-1i}_{\enspace 2j}
=(R^{-1})^{1i}_{2j} = 0
\end{array}
\end{equation}
and so on for their pairing with products of the $\ct^{i}_{j}$. Therefore,
setting these elements to zero yields a quotient bialgebra $U(R)$ of
$\tilde{U}(R)$
\begin{equation}
\cl^{+}=\begin{pmatrix}
\cl^{+0}_{\enspace 0} & 0 & 0\\
0 & \cl^{+1}_{\enspace 1} & \cl^{+1}_{\enspace 2}\\
0 & 0 & \cl^{+2}_{\enspace 2}
\end{pmatrix}, \qquad
\cl^{-}=\begin{pmatrix}
\cl^{-0}_{\enspace 0} & 0 & 0\\
0 & \cl^{-1}_{\enspace 1} & 0\\
0 & \cl^{-2}_{\enspace 1} & \cl^{-2}_{\enspace 2}
\end{pmatrix}    
\end{equation}   
Therefore, the initial pairing $\langle A(R),\tilde{U}(R)\rangle$
descends to $\langle \ca_{r,s}, U(R)\rangle$. So, for $\cu_{r,s}$ (or
$U(R)$) we make the following ansatz for the $\cl^{\pm}$
matrices:

\begin{eqnarray*}
\cl^{+} &=& c^{+}r \begin{pmatrix}
s^{-\frac{1}{2} (\tf-H_{2}-1)}r^{\frac{1}{2} (\tf-H_{1}-1)} & 0 & 0\\   
0 & s^{-\frac{1}{2}(\tf-H_{1}+1)}r^{\frac{1}{2}(-\tf+H_{2}-1)}
& r^{-1}\lambda \tc\\
0 & 0 & s^{-\frac{1}{2}(\tf+H_{1}-1)}r^{\frac{1}{2}(-\tf-H_{2}-1)}\\
\end{pmatrix}
\\
\cl^{-} &=& c^{-}r^{-1} \begin{pmatrix}
s^{-\frac{1}{2} (\tf-H_{2}-1)}r^{-\frac{1}{2} (\tf-H_{1}-1)} & 0 & 0\\
0 & s^{-\frac{1}{2}(\tf-H_{1}+1)}r^{-\frac{1}{2}(-\tf+H_{2}-1)} & 0\\
0 & -r \lambda \tb & 
s^{-\frac{1}{2}(\tf+H_{1}-1)}r^{-\frac{1}{2}(-\tf-H_{2}-1)}\\
\end{pmatrix}  
\end{eqnarray*}

where $H_{1}=\ta+\td$, $H_{2}=\ta-\td$, and $\{ \ta, \tb, \tc, \td, \tf
\}$
is the set of generating elements of the dual algebra. This is consistent
with the action on the generators of $\ca_{r,s}$ and gives the
correct duality pairings. More conveniently,

\begin{equation}
\cl^{+}= \begin{pmatrix}
J & 0 & 0\\
0 & M & P\\
0 & 0 & N\\
\end{pmatrix} \quad \text{and} \quad
\cl^{-}= \begin{pmatrix}
J' & 0 & 0\\
0 & M' & 0\\
0 & Q & N'\\
\end{pmatrix}
\end{equation}
where
\begin{equation}
\begin{array}{ll}
&J=s^{-\frac{1}{2} (\tf-H_{2}-1)}r^{\frac{1}{2} (\tf-H_{1}+1)}\\
&M=s^{-\frac{1}{2} (\tf-H_{1}+1)}r^{\frac{1}{2} (-\tf+H_{2}+1)}\\
&N=s^{-\frac{1}{2}(\tf+H_{1}-1)}r^{\frac{1}{2} (-\tf-H_{2}+1)}\\
&J'=s^{-\frac{1}{2} (\tf-H_{2}-1)}r^{-\frac{1}{2} (\tf-H_{1}+1)}\\
&M'=s^{-\frac{1}{2} (\tf-H_{1}+1)}r^{-\frac{1}{2} (-\tf+H_{2}+1)}\\
&N'=s^{-\frac{1}{2}(\tf+H_{1}-1)}r^{-\frac{1}{2} (-\tf-H_{2}+1)}
\end{array}
\end{equation}
and
\begin{equation}
\begin{array}{ll}
&P=\lambda \tc\\
&Q=-\lambda \tb
\end{array}  
\end{equation} 
These can also be arranged in terms of smaller $L^{+}$ and $L^{-}$
matrices
\begin{equation}
\begin{array}{llll}
&\cl^{+}=c^{+}\begin{pmatrix}
J & 0\\
0 & L^{+}
\end{pmatrix} \quad &\text{where} \quad
&L^{+}=\begin{pmatrix}
M & P\\
0 & N
\end{pmatrix}\\
&\cl^{-}=c^{-}\begin{pmatrix}
J' & 0\\
0 & L^{-}
\end{pmatrix} \quad &\text{where} \quad
&L^{-}=\begin{pmatrix}
M' & 0\\
Q & N'  
\end{pmatrix}
\end{array}
\end{equation}

\subsection*{Commutation relations of the dual}
The dual algebra is generated by $\cl^{\pm}$ functionals which satisfy the
$q$-commutation relations (the so-called $R\cl\cl$ relations)
\begin{eqnarray}
R_{12}\cl^{\pm}_{2}\cl^{\pm}_{1} &=& \cl^{\pm}_{1}\cl^{\pm}_{2}R_{12}\\
R_{12}\cl^{+}_{2}\cl^{-}_{1} &=& \cl^{-}_{1}\cl^{+}_{2}R_{12}
\end{eqnarray}
where $\cl^{\pm}_{1}=\cl^{\pm}\ot \ident$ and  $\cl^{\pm}_{2}=\ident\ot
\cl^{\pm}$. Since $\ca_{r,s}$ is a quotient Hopf algebra, it
is necessary to amend the $R$-matrix to eliminate relations that are
inconsistent with the quotient structure. Consequently, the $R$-matrix for
the $R\cl\cl$ relations is different from the one used in the $R\ct\ct$
relations. The $R\cl\cl$ relations are constructed with the $R$-matrix:
\begin{equation}
R_{12}=c^{-}{\langle \cl^{-},\ct \rangle}^{-1}
=\begin{pmatrix}
r & 0 & 0 & 0\\
0 & \mathsf{S}^{-1} & 0 & 0\\
0 & 0 & \mathsf{S} & 0\\
0 & 0 & 0 & R_{r}
\end{pmatrix}
\end{equation}
Evaluating $\cl^{\pm}_{1}$, $\cl^{\pm}_{2}$ matrices and substituting in
the above $R\cl\cl$- relations yields the dual algebra commutation 
relations.
From
$R_{12}\cl^{-}_{2}\cl^{-}_{1} = \cl^{-}_{1}\cl^{-}_{2}R_{12}$ and
$R_{12}\cl^{+}_{2}\cl^{+}_{1} = \cl^{+}_{1}\cl^{+}_{2}R_{12}$
we obtain
\begin{eqnarray}
R_{r}L^{-}_{2}L^{-}_{1} &=& L^{-}_{1}L^{-}_{2}R_{r}\\
R_{r}L^{+}_{2}L^{+}_{1} &=& L^{+}_{1}L^{+}_{2}R_{r}
\end{eqnarray}
\begin{equation}
\begin{array}{lll}
& MJ=JM \qquad & M'J'=J'M'\\
& NJ=JN \qquad & N'J'=J'N'\\
& PJ=sJP \qquad & J'Q=sQJ'
\end{array}
\end{equation}
where
\begin{equation}
\begin{array}{llllll}
&R_{r}L^{-}_{2}L^{-}_{1}=L^{-}_{1}L^{-}_{2}R_{r} &\Longrightarrow
&QM'=rM'Q,\quad N'Q=rQN' &\text{and} &N'M'=M'N'\\
&R_{r}L^{+}_{2}L^{+}_{1}=L^{+}_{1}L^{+}_{2}R_{r} &\Longrightarrow
&PM=rMP,\quad NP=rPN &\text{and} &NM=MN
\end{array}
\end{equation}
In addition, the cross relation 
$R_{12}\cl^{+}_{2}\cl^{-}_{1} = \cl^{-}_{1}\cl^{+}_{2}R_{12}$ yields
\begin{equation}
\begin{array}{llll}
& NJ'=J'N \qquad & MJ'=J'M \qquad & PJ'=sJ'P\\
& N'J=JN' \qquad & M'J=JM' \qquad & JQ=sQJ
\end{array}
\end{equation}
and $R_{r}L^{+}_{2}L^{-}_{1} = L^{-}_{1}L^{+}_{2}R_{r}$ which further
implies
\begin{equation}
QP-PQ=-\lambda(N'M-NM')
\end{equation}
Simplifying the above, we get the following commutation relations
\begin{equation}
\begin{array}{lll}
&[\ta,\tb]=\tb, \quad &[\ta,\tc]=-\tc\\
&[\td,\tb]=-\tb, \quad &[\td,\tc]=\tc\\
&[\ta,\td]=0, \quad &[\tf, \bullet]=0\\
\end{array}
\end{equation}
and 
\begin{equation}
[\tb,\tc]=
\frac{r^{\ta-\td}s^{-\tf}-{r^{-(\ta-\td)}s^{-\tf}}}{r-r^{-1}}= 
\frac{r^{\gamma \tf}}{r-r^{-1}}[r^{\ta-\td}-r^{-(\ta-\td)}]
\end{equation}
where $\gamma =\frac{\ln s}{\ln r}$. So, we obtain a single-parameter
deformation of $U(gl(2))\ot U(u(1))$ as an algebra. Including the
coproduct, we again obtain a semidirect product
$U_{r}(gl(2))\underset{s}{\rtimes} U(u(1))$, as expected.

\section{Constructive calculus}

In order to investigate the differential geometric structure of the
$(r,s)$-deformation, $\ca_{r,s}$, of $GL(2)\ot GL(1)$, we use Jur\u{c}o's
constructive procedure \cite{jurco} based on the $R$-matrix
formulation. This method has so far been applied only to full matrix
quantum groups but we demonstrate here that it works equally well for
appropriate quotients of these. For $\ca_{r,s}$, we obtain a first order
bicovariant differential calculus employing the ansatz for $\cl^{\pm}$
introduced in Sec. III.

\subsection{One-forms}

Let $\{ \omega \}$ be the basis of all left-invariant quantum 
one-forms. So, we have
\begin{equation}
\cop_{L}(\omega)= \ident\ot \omega
\end{equation}
This defines the left action on the bimodule $\Gamma$ (space of quantum
one-forms). The bimodule $\Gamma$ is further characterised by the
commutation relations between $\omega$ and $a\in \ca$ ($\equiv\ca_{r,s}$),
\begin{equation}
\omega a = (f\ast a) \omega
\end{equation}
The left convolution product is
\begin{equation}
f\ast a = (\ident\ot f) \cop (a)
\end{equation}
where $f\in \ca'(=\text{Hom}(\ca, \mathbf{C} ))$ belongs to the dual. This
means
\begin{equation}
\omega a = (\ident\ot f) \cop (a) \omega
\end{equation}
Now the linear functional $f$ is defined in terms of the $\cl^{\pm}$
matrices as
\begin{equation}
f = S(\cl^{+})\cl^{-}
\end{equation}
Thus we have
\begin{equation}
\omega a = [(\ident\ot S(\cl^{+})\cl^{-}) \cop (a)] \omega
\end{equation}
In terms of components,
\begin{equation}
\omega_{ij} a = [(\ident\ot S(l^{+}_{ki})l^{-}_{jl}) \cop (a)] \omega_{kl}
\end{equation}
using the expressions $\cl^{\pm} = l^{\pm}_{ij}$ and $\omega =
\omega_{ij}$ where $i, j = 1 .. 3$. For $\Gamma$ to be a bicovariant
bimodule, the right coaction is given by
\begin{equation}
\cop_{R} (\omega) = \omega\ot M
\end{equation}
where functionals $M$ are defined in terms of the matrix of generators
$\ct$,
\begin{equation}
M = \ct S(\ct)
\end{equation}
Again, in component form, we can write
\begin{equation}
\cop_{R} (\omega_{ij}) = \omega_{kl}\ot t_{ki}S(t_{jl})
\end{equation}
Using the above formulas, we obtain the commutation relations of all the
left-invariant one forms with the generating elements $\{ a,b,c,d,f \}$ of
$\ca_{r,s}$:
\begin{equation}
\begin{array}{ll}
\omega^{0} a = a\omega^{0} & \qquad \qquad
\omega^{0} b = b\omega^{0}\\
\omega^{1} a = r^{-2}a\omega^{1} & \qquad \qquad
\omega^{1} b = b\omega^{1}\\
\omega^{+} a = r^{-1}a\omega^{+} & \qquad \qquad
\omega^{+} b = r^{-1}b\omega^{+} - \lambda r^{-1}a\omega^{1}\\
\omega^{-} a = r^{-1}a\omega^{-} -\lambda r^{-1}b\omega^{1} & \qquad
\qquad
\omega^{-} b = r^{-1}b\omega^{-}\\
\omega^{2} a = a\omega^{2} -\lambda b\omega^{+} & \qquad \qquad
\omega^{2} b = r^{-2}b\omega^{2} -\lambda r^{-1}a\omega^{-} +
\lambda^{2} b\omega^{1}\\
\end{array}   
\end{equation}
\begin{equation}
\begin{array}{ll}
\omega^{0} c = c\omega^{0} & \qquad \qquad
\omega^{0} d = d\omega^{0}\\
\omega^{1} c = r^{-2}c\omega^{1} & \qquad \qquad
\omega^{1} d = d\omega^{1}\\
\omega^{+} c = r^{-1}c\omega^{+} & \qquad \qquad
\omega^{+} d = r^{-1}d\omega^{+} - \lambda r^{-1}c\omega^{1}\\
\omega^{-} c = r^{-1}c\omega^{-} -\lambda r^{-1}d\omega^{1} & \qquad
\qquad
\omega^{-} d = r^{-1}d\omega^{-}\\
\omega^{2} c = c\omega^{2} -\lambda d\omega^{+} & \qquad \qquad
\omega^{2} d = r^{-2}d\omega^{2} -\lambda r^{-1}c\omega^{-} +
\lambda^{2} d\omega^{1}
\end{array}
\end{equation}
\begin{equation}
\begin{array}{l}
\omega^{0} f = r^{-2}f\omega^{0}\\
\omega^{1} f = f\omega^{1}\\
\omega^{+} f = sf\omega^{+}\\
\omega^{-} f = s^{-1}f\omega^{-}\\
\omega^{2} f = f\omega^{2}
\end{array}
\end{equation}
where $\omega^{0}=\omega_{11}, \omega^{1}=\omega_{22},
\omega^{+}=\omega_{23}, \omega^{-}=\omega_{32}, \omega^{2}=\omega_{33}$
and the components $\omega_{12}, \omega_{13}, \omega_{21}, \omega_{31}$
have null contribution, given the structure of the $\ct$ matrix (i.e.,
$t_{12}=t_{13}=t_{21}=t_{31}=0$).

\subsection{Vector fields}

The linear space $\Gamma$ (space of all left invariant
one-forms) contains a bi-invariant element $\tau = \sum_{i}\omega_{ii}$
which can be used to define a derivative on $\ca$. For $a\in \ca$, one
sets
\begin{equation}
\mathbf{d} a = \tau a - a \tau
\label{tau}
\end{equation}
Now
\begin{equation}
\omega_{ii} a = [(\ident\ot S(l^{+}_{ki})l^{-}_{il}) \cop (a)] \omega_{kl}
\end{equation}
So
\begin{equation}
\mathbf{d}a = [(\ident\ot \chi_{kl}) \cop (a)] \omega_{kl}
\end{equation}
where $\chi_{kl} = S(l^{+}_{ki})l^{-}_{il} - \delta_{kl}\cnt$, $\cnt$
being the counit. Denote
\begin{equation}
\chi_{ij} = S(l^{+}_{ik})l^{-}_{kj} - \delta_{ij}\cnt
\end{equation}
or more compactly
\begin{equation}
\chi = S(\cl^{+})\cl^{-} - \ident\cnt
\end{equation}
the matrix of left-invariant vector fields $\chi_{ij}$ on $\ca$. The
action of the vector fields on the generating elements is
\begin{eqnarray}
\chi_{ij} a &=& (S(l^{+}_{ik})l^{-}_{kj} - \delta_{ij}\cnt)a\\
\chi_{ij} a &=& \langle S(l^{+}_{ik})l^{-}_{kj}, a\rangle -
\delta_{ij}\cnt (a)
\end{eqnarray}
Explicitly, we obtain
\begin{equation}
\begin{array}{ll}
\chi_{0}(a) = 0 & \qquad \chi_{0}(b) = 0\\
\chi_{1}(a) = r^{-2}-1 & \qquad \chi_{1}(b) = 0\\
\chi_{+}(a) = 0 & \qquad \chi_{+}(b) = 0\\
\chi_{-}(a) = 0 & \qquad \chi_{-}(b) = -(r-r^{-1})\\
\chi_{2}(a) = 0 & \qquad \chi_{2}(b) = 0\\
\end{array}
\end{equation}
\begin{equation}
\begin{array}{ll}
\chi_{0}(c) = 0 & \qquad \chi_{0}(d) = 0\\
\chi_{1}(c) = 0 & \qquad \chi_{1}(d) = (r-r^{-1})^{2}\\
\chi_{+}(c) = -(r-r^{-1}) & \qquad \chi_{+}(d) = 0\\
\chi_{-}(c) = 0 & \qquad \chi_{-}(d) = 0\\
\chi_{2}(c) = 0 & \qquad \chi_{2}(d) = r^{-2}-1
\end{array}    
\end{equation}
\begin{equation}
\begin{array}{l}
\chi_{0}(f) = r^{-2}-1\\
\chi_{1}(f) = 0\\
\chi_{+}(f) = 0\\
\chi_{-}(f) = 0\\
\chi_{2}(f) = 0
\end{array}   
\end{equation}

where $\chi_{0}=\chi_{11}, \chi_{1}=\chi_{22}, \chi_{+}=\chi_{23},
\chi_{-}=\chi_{32}, \chi_{2}=\chi_{33}$ and again (by previous argument)
the components $\chi_{12}, \chi_{13}, \chi_{21}, \chi_{31}$ have null
contribution. The left convolution products are given as

\begin{equation}
\begin{array}{ll}
\chi_{0} \ast a = 0 & \qquad \chi_{0} \ast b = 0\\
\chi_{1} \ast a = (r^{-2}-1)a & \qquad \chi_{1} \ast b =
((r-r^{-1})^{2})b\\
\chi_{+} \ast a = -(r-r^{-1})b & \qquad \chi_{+} \ast b = 0\\
\chi_{-} \ast a = 0 & \qquad \chi_{-} \ast b = -(r-r^{-1})a\\
\chi_{2} \ast a = 0 & \qquad \chi_{2} \ast b = (r^{-2}-1)b\\
\end{array}
\end{equation}   
\begin{equation}
\begin{array}{ll}
\chi_{0} \ast c = 0 & \qquad \chi_{0} \ast d = 0\\
\chi_{1} \ast c = (r^{-2}-1)c & \qquad \chi_{1} \ast d =
((r-r^{-1})^{2})d\\
\chi_{+} \ast c = -(r-r^{-1})d & \qquad \chi_{+} \ast d = 0\\
\chi_{-} \ast c = 0 & \qquad \chi_{-} \ast d = -(r-r^{-1})c\\
\chi_{2} \ast c = 0 & \qquad \chi_{2} \ast d = (r^{-2}-1)d
\end{array}   
\end{equation}
\begin{equation}
\begin{array}{l}
\chi_{0} \ast f = (r^{-2}-1)f\\
\chi_{1} \ast f = 0\\
\chi_{+} \ast f = 0\\
\chi_{-} \ast f = 0\\
\chi_{2} \ast f = 0
\end{array}   
\end{equation}

\subsection{Exterior derivatives}

Using $\mathbf{d} a =\sum_{i}(\chi_{i} \ast a) \omega^{i}$ for
$a \in \ca$, we obtain the action of the exterior derivatives:
\begin{eqnarray}
\mathbf{d} a &=& (r^{-2}-1)a\omega^{1} - \lambda
b\omega^{+}\\
\mathbf{d} b &=& \lambda^{2}b\omega^{1} - \lambda
a\omega^{-} + (r^{-2}-1)b\omega^{2}\\
\mathbf{d} c &=& (r^{-2}-1)c\omega^{1} - \lambda
d\omega^{+}\\
\mathbf{d} d &=& \lambda^{2}d\omega^{1} - \lambda  
c\omega^{-} + (r^{-2}-1)d\omega^{2}\\
\mathbf{d} f &=& (r^{-2}-1)f\omega^{0}
\end{eqnarray}
where $\lambda=r-r^{-1}$. The exterior derivative $\mathbf{d}: \ca \lra
\Gamma$ satisfies the  Leibniz rule and $\mathbf{d}\ca$ generates $\Gamma$
as a left $\ca$-module. This then defines a first-order differential
calclulus $(\Gamma, \mathbf{d})$ on $\ca_{r,s}$. Furthermore, the calculus
is bicovariant due to the coexistence of the left and the right actions
\begin{eqnarray}
&& \cop_{L}:\Gamma \lra \ca \ot \Gamma\\
&& \cop_{R}:\Gamma \lra \Gamma \ot \ca
\end{eqnarray}
since $\mathbf{d}$ has the invariance property
\begin{eqnarray}
&& \cop_{L}\mathbf{d}=(\ident \ot \mathbf{d})\cop\\
&& \cop_{R}\mathbf{d}=(\mathbf{d} \ot \ident)\cop       
\end{eqnarray}
The bicovariance holds also due to the existence of the bi-invariant
element $\tau = \sum_{i}\omega_{ii}$ (eqn.(\ref{tau})) of the linear space
of left-invariant one-forms. If we rewrite the derivatives $\{
\mathbf{d}a,\mathbf{d}b,\mathbf{d}c,\mathbf{d}d,\mathbf{d}f \}$ as $\{ 
\mathbf{d}(f^{N}a),\mathbf{d}(f^{N}b),\mathbf{d}(f^{N}c),\mathbf{d}(f^{N}d) 
\}$, i.e., reducing from the five-dimensional to the four-dimensional
algebra, then the latter set of exterior derivatives provides a
realisation of the differential calculus on the biparametric
$(p,q)$-deformation of $GL(2)$, i.e., $\GLpqtwo$, with the defining
relations between the two sets of deformation parameters $(p,q)$ and
$(r,s)$ as before. Furthermore, the differential calculus also respects
the cross-product structure of $\ca_{rs}$. It can be checked (using the
Leibniz rule) that
\begin{equation}
\mathbf{d}(af-fa)=0, \quad \mathbf{d}(cf-sfc)=0, \quad
\mathbf{d}(bf-s^{-1}fb)=0, \quad \mathbf{d}(df-fd)=0, \quad
\end{equation}
which is consistent with the cross relations (\ref{grsf}).

\section{Jordanian analogue}

It was shown in \cite{deeps} that the $\ca_{r,s}$ deformation could be
contracted (by means of singular limit of similarity transformations) to
obtain a nonstandard or Jordanian analogue, say $\ca_{m,k}$, with
deformation parameters $\{ m,k \}$ and the associated $R$-matrix is
triangular. In analogy with $\ca_{r,s}$, $\ca_{m,k}$ can also be
considered as the semidirect or cross-product
$GL_m(2)\underset{k}{\rtimes} \mathbb{C}[f,f^{-1}]$ where
$GL_m(2)=\mathbb{C}[a,b,c,d]$ modulo the relations 
\begin{equation}
\begin{array}{ll}
{[}c,d{]} = -mc^{2}, & \qquad [c,b] = -m(ac+cd) = -m(ca+dc) \\
{[}c,a{]} = -mc^{2}, & \qquad [d,a] = -m(d-a)c = -mc(d-a)
\end{array}
\label{gmkabcd1}
\end{equation}
\begin{equation}
\begin{array}{l}
{[}d,b{]} = -m(d^{2}-\delta) \\
{[}b,a{]} = -m(\delta-a^{2})
\end{array}  
\label{gmkabcd2}
\end{equation}
where $\delta=ad-bc+mac=ad-cb-mcd$, and $\mathbb{C}[f,f^{-1}]$ has the
cross relations 
\begin{equation}
\begin{array}{ll}
{[}f,a{]} = kcf, & \qquad[f,b] = k(df-fa) \\
{[}f,c{]} = 0, & \qquad[f,d] = -kcf
\end{array}  
\label{gmkf}
\end{equation}
Thus, $\ca_{m,k}
\simeq GL_m(2)\underset{k}{\rtimes} \mathbb{C}[f,f^{-1}]$ can also 
be interpreted as a skew Laurent polynomial ring $GL_{m}[f,f^{-1};\sigma]$
where $\sigma$ is the automorphism given by the action of element $f$ on
$GL_m(2)$. The (left) action is given by
\begin{equation}
f\lact a=a+kc, \qquad f\lact b=b+k(d-a)-k^{2}c, \qquad f\lact c=c, \qquad
f\lact d=d-kc
\end{equation}

\section{Discussion}

In this article, we have investigated the algebro-geometric structure of
the
biparametric quantum deformation of $GL(2)\ot GL(1)$, namely, 
$\ca_{r,s}$. A particular feature  of this deformation is that it has an
interpretation as a semidirect or cross-product algebra. We exhibit this
cross-product structure and establish a picture of duality in this
setting. Using the $R$-matrix formalism, we have given an explicit
derivation of the corresponding dual agebra, i.e., the quantised
universal enveloping algebra and also constructed a bicovariant
differential calculus. The dual algebra obtained via $R$-matrices is
isomorphic to the dual algebra obtained by the cross-product construction. 
We note that the differential calculus satisfies the required axioms,
contains the calculus on $\GLqtwo$ and our results match with those given
in \cite{cast}. Besides, the calculus is also consistent with the
cross-product structure of $\ca_{r,s}$. We expect that the calculus could as
well be obtained by projection from the calculus on multiparameter
$q$-deformed $GL(3)$.  The differential calculus obtained on $\ca_{r,s}$
enables us to investigate the associated gauge theory from a
noncommutative perspective. It would be useful to repeat the analysis
presented in this paper for the biparametric Jordanian deformation of
$GL(2)\ot GL(1)$ obtained in \cite{deeps} and also to investigate
corresponding hybrid $(q,h)$-deformations \cite{bhp, hybrid}. Furthermore,
it would indeed be interesting to generalise the setting to the case of
coloured quantum and Jordanian deformations.

\section*{Acknowledgements}

I am grateful to Profs. Shahn Majid and Konrad Schm\"udgen for several
fruitful discussions.

\end{document}